\documentclass[a4paper,12pt]{amsart}
\usepackage[labelsep=colon]{caption}
\usepackage{float}
\restylefloat{figure}
\usepackage{graphicx}
\usepackage[colorlinks=true]{hyperref}
\hypersetup{allcolors=[rgb]{0.9,0.1,0.4}}
\usepackage{mathrsfs}
\usepackage[sort&compress,numbers]{natbib}
\usepackage{slashbox}
\usepackage[hang]{subfigure}
\usepackage{tikz}
\usetikzlibrary{arrows,%
  decorations.markings,%
  decorations.pathmorphing,%
  matrix,%
  shapes}
\usepackage{wrapfig}

\hypersetup{
  pdftitle={Topology and Higher Concurrencies},%
  pdfauthor={Nils A. Baas},%
}

\tikzstyle{midarrow} = [decoration={markings,mark=at position .5 with
  \arrow{>}}, postaction={decorate}]

\DeclareMathOperator{\source}{source}
\DeclareMathOperator{\supp}{supp}
\DeclareMathOperator{\target}{target}
\newcommand{\boks}{\mathop{}\!\square\mathop{}\!}
\newcommand{\C}{\mathscr{C}}

\renewcommand{\H}{\mathscr{H}}
\newcommand{\hc}{\hat{c}}

\renewcommand{\P}{\mathscr{P}}

\title{Topology and Higher Concurrencies}
\author[N.A. Baas]{Nils A.\ Baas}
\date{May 3, 2018}
\address{Department of Mathematical Sciences, NTNU, NO-7491 Trondheim,
  Norway}
\email{nils.baas@ntnu.no}

\begin{document}
\maketitle

In general concurrence means the combined action of agents or causes,
often in the form of simultaneous occurrence. It is also a form of
correlation. Examples are cofiring of neurons, coexpression of genes,
events happening at the same time, etc. In this paper we want to
extend the concept of concurrency to a higher version fitting into the
general scheme of higher structures (hyperstructures) as described and
studied in \cite{Hnano, Cog, HAM, HTD, NSCS, NS, SO, HOA, HOS, MHS,
  PHS, BS, BF, BSS}. For example neurons may fire simultaneously in
groups, even groups of groups, or it could be brain domains active in
various tasks being performed and measured by fMRI techniques.

We are inspired by th work of Ellis and Klein \cite{EK}, and the
framework for neural codes by Curto and Itskov \cite{Curto}.

The idea of introducing higher concurrencies in general applies to all
kinds of events. We will here illustrate the idea in the framework of
neural codes as in \cite{Curto} and introduce higher neural codes. We
will show how the use of topology extends to this case as well.

By a neural code word of length $n$ we mean a sequence
\begin{equation*}
  C = (\hc_1,\ldots,\hc_n), \quad \hc_i \in \{0,1\}.
\end{equation*}
A code is a finite collection of such words
\begin{equation*}
  \C = \{c_i\}.
\end{equation*}
We follow the notation of \cite{Curto}.\\

Hence $\C \subseteq 2^{[n]} = \P[n]$ or $\C \in \P^2 [n]$.\\

We define the support as follows:
\begin{align*}
  \supp C &= \{i \in [n] \mid \hc_i = 1\} \subseteq [n]\\
  \supp \C &= \{\supp(C) \mid C\in \C\} \subseteq 2^n = \P[n]\\
  \supp \C &\in \P^2[n].
\end{align*}
Such words come from neural data (spike trains) as follows:
\begin{center}
  \begin{tabular}{c|cc}
    \backslashbox{neurons}{time} & $\cdots$ & t\\ \hline
    $N_1$ && $c_1$\\
    $N_2$ && $c_2$\\
    $\vdots$ && $\vdots$\\
    $N_k$ && $c_n$
  \end{tabular}
\end{center}
We are naturally interested in cofiring patterns, hence we pay
attention to $\supp \C$.

We are looking for (topological) structure in the coding data.  Let us
form:\\
\begin{equation*}
  \Delta(\C) = \{\sigma \subseteq [n] \mid \sigma \subseteq \supp(C)
  \text{ for some $C$}\}
\end{equation*}
\mbox{}\\
\noindent --- this is the smallest simplicial complex that contains
$\supp \C$.

The question we could like to address is: \emph{Is there some higher
  order (topological) information in the code data, how to extract
  it?}

All thinking about higher order structures in set theory, type theory,
logic, higher categories, etc.\ is in some way often a sophisticated
reflection of the basis hierarchy:
\begin{equation*}
  X, \P(X), \P^2(X),\ldots,\P^n(X),\ldots
\end{equation*}
Hyperstructures, see \cite{Hnano, Cog, HAM, HTD, NSCS, NS, SO, HOA,
  HOS, MHS, PHS, BS, BF, BSS} for references, represent a highly
refined version of this.  Let us pursue this in our coding context
combining with the ideas from hyperstructures, hence constructing an
interesting example.

Put $\C_1 = \supp \C\in \P^2[n]$, an element represents a cofiring or a
concurency
\begin{equation*}
  (101101110\ldots 0) \mapsto (111111)
\end{equation*}
pattern.  Then we may look for patterns of patterns.

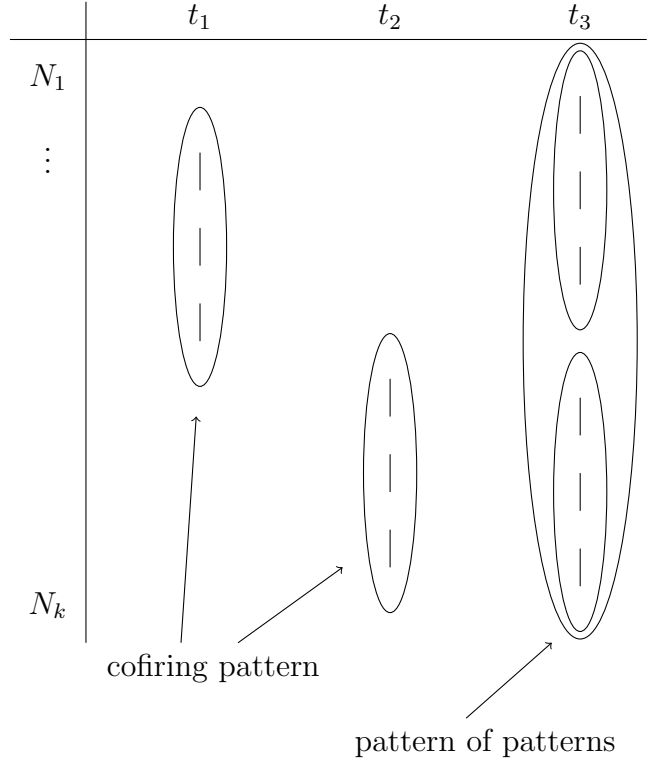
\begin{center}
  \begin{tikzpicture}
    \draw (0,0) -- (8.375,0);
    \draw (1,0.5) -- (1,-8);
    \node at (0.5,-0.5){$N_1$};
    \node at (0.5,-1.5){$\vdots$};
    \node at (0.5,-7.5){$N_k$};
    \node[above] at (2.5,0){$t_1$};
    \node[above] at (5,0){$t_2$};
    \node[above] at (7.5,0){$t_3$};

    \begin{scope}[xshift=2.5cm,yshift=-1.5cm]
      \draw (0,0) -- (0,-0.5);
      \draw (0,-1) -- (0,-1.5);
      \draw (0,-2) -- (0,-2.5);
      \draw (0,-1.25) ellipse (0.35cm and 1.85cm);
    \end{scope}

    \begin{scope}[xshift=5cm,yshift=-4.5cm]
      \draw (0,0) -- (0,-0.5);
      \draw (0,-1) -- (0,-1.5);
      \draw (0,-2) -- (0,-2.5);
      \draw (0,-1.25) ellipse (0.35cm and 1.85cm);
    \end{scope}

    \begin{scope}[xshift=7.5cm,yshift=-0.75cm]
      \draw (0,0) -- (0,-0.5);
      \draw (0,-1) -- (0,-1.5);
      \draw (0,-2) -- (0,-2.5);
      \draw (0,-1.25) ellipse (0.35cm and 1.85cm);
    \end{scope}

    \begin{scope}[xshift=7.5cm,yshift=-4.75cm]
      \draw (0,0) -- (0,-0.5);
      \draw (0,-1) -- (0,-1.5);
      \draw (0,-2) -- (0,-2.5);
      \draw (0,-1.25) ellipse (0.35cm and 1.85cm);
    \end{scope}

    \begin{scope}[xshift=7.5cm, yshift=-4cm]
      \draw (0,0) ellipse (0.75cm and 3.95cm);
    \end{scope}

    \draw[->] (2.25,-8) -- (2.45,-5);
    \draw[->] (3,-8) -- (4.375,-7);
    \node[below] at (2.675,-8){cofiring pattern};

    \draw[->] (6,-9) -- (7.125,-8);
    \node[below] at (6.25,-9){pattern of patterns};
  \end{tikzpicture}
\end{center}

To be more precise a pattern exists when $\supp C\in \C_1$ is realized
at a certain time.

If several patterns are realized or occur at another time ($t_3$) we say
that they form a second order pattern.  In the language of
hyperstructures
\begin{itemize}
\item[] \emph{Observer 1:} individual neurons firing (or other
  activity)\\
\item[] \emph{Bond 1:} cofiring of neurons (coactivity)\\
\item[] \emph{Observer 2:} firing of a pattern (activity of
  patterns)\\
\item[] \emph{Bond 2:} cofiring of patterns (coactivity of patterns).
\end{itemize}
This clearly continues, and should be managable combinatorially due to
a reasonably limited number of cofiring patterns.

Hence we put
\begin{equation*}
  \C_2 = \supp \C_1 \in \P^3 [n]
\end{equation*}
representing cofiring patterns
\begin{center}
  \begin{tikzpicture}[scale=0.5]
    \draw (0,0) -- (0,-0.5);
    \draw (0,-1) -- (0,-1.5);
    \draw (0,-2) -- (0,-2.5);
    \draw (0,-1.25) ellipse (0.35cm and 1.85cm);
  \end{tikzpicture}
\end{center}
\begin{equation*}
  \C_3 = \supp \C_2 \in \P^4 [n]
\end{equation*}
representing
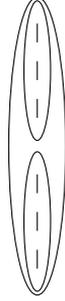
\begin{center}
  \begin{tikzpicture}[scale=0.5]
    \begin{scope}[xshift=0cm,yshift=-0.75cm]
      \draw (0,0) -- (0,-0.5);
      \draw (0,-1) -- (0,-1.5);
      \draw (0,-2) -- (0,-2.5);
      \draw (0,-1.25) ellipse (0.35cm and 1.85cm);
    \end{scope}

    \begin{scope}[xshift=0cm,yshift=-4.75cm]
      \draw (0,0) -- (0,-0.5);
      \draw (0,-1) -- (0,-1.5);
      \draw (0,-2) -- (0,-2.5);
      \draw (0,-1.25) ellipse (0.35cm and 1.85cm);
    \end{scope}

    \begin{scope}[xshift=0cm, yshift=-4cm]
      \draw (0,0) ellipse (0.75cm and 3.95cm);
    \end{scope}
  \end{tikzpicture}
\end{center}
Iteratively:\\
\begin{equation*}
  \C_k = \supp \C_{k - 1} \in \P^{k + 1} [n].
\end{equation*}
\mbox{}\\

If we include $k$ levels of structure then
\begin{equation*}
  \C = \{\C_1,\C_2,\ldots,\C_k\}
\end{equation*}
form a hyperstructure with ``boundary'' maps
\begin{equation*}
  \partial_i \colon \C_{i + 1} \to \P(\C_i), \quad \partial_i (c_{i +
    1}) = \{c_i^j\} \quad [c_{i + 1} = \{c_i^1,\ldots,c_i^{i_k}\}]
\end{equation*}
dissolving the bond (as for geometric cobordisms!).  Let us call the
hyperstructure $\H(\C)$.

The hyperstructure is a useful way to organize the code data.  How to
extract topological information out of this?

We may associate the clique topology complex to each level
\begin{equation*}
  \Delta(\C) \colon \Delta (\C_k) \xleftarrow{\delta_{k - 1}} \Delta
  (\C_{k - 1}) \leftarrow \cdots \xleftarrow{\delta_1} \Delta (\C_1)
\end{equation*}
where the $\delta_i$'s are induced from the $\partial_i$'s. We assume
these are well-defined.

Clearly this sequence gives more information than just $\Delta
(\C_1)$.  For various filtrations (frequency, ...) we may apply
persistent homology --- $PH$ --- to the sequence and get:
\begin{equation*}
  PH (\Delta(\C_k)) \xleftarrow{\delta_{k - 1}^\ast} PH(\Delta(\C_{k
    - 1})) \leftarrow \cdots \xleftarrow{\delta_1^\ast} PH(\Delta(\C_1))
\end{equation*}
giving a sequence of barcodes (persistent diagrams) and maps as an
invariant.

But it does not capture the fact that many cofiring patterns are not
simplicial, in the sense that if a collection cofires, not necessarily
all subcollections will cofire (at some time).  This is
organizationally captured in the hyperstructure.

How can this get reflected topologically?

We will associate a space --- a simplicial complex, the nerve --- to
the hyperstructure and then study its shape via (persistent)
homology.  Let us return to the hyperstructure of coding (firing)
patterns:
\begin{equation*}
  \H(\C) = \{\C_1,\C_2,\ldots,\C_k\}.
\end{equation*}
We will construct ``the nerve of $\H(\C)$''.  What should the
simplexes be?  

Let us use the construction of the classifying space of a category as
a model, where the $p$-simplexes are compositions of morphisms of
length $p$.  We have to define composition of bonds at any level.  We
define compositions of bonds by ``gluing'' similar to morphisms in
categories and higher categories.
\begin{equation*}
  \bullet \xrightarrow{f} \bullet \xrightarrow{g} \bullet
\end{equation*}
$g\circ f$ is defined if $\target f = \source g$.  For $n$-morphisms
we may have that at level $p$:
\begin{equation*}
  \target_p^n (f) = \source_p^n (g)
\end{equation*}
``gluing'' at level $p$, giving a composition
\begin{equation*}
  f \boks_p^n\; g
\end{equation*}

For a hyperstructure we search for compatibility conditions in order
to glue.  If we have two $n$-bonds:
\begin{equation*}
  a_n \text{ and } b_n
\end{equation*}
and
\begin{equation*}
  \partial_p \circ \cdots \circ \partial_{n - 1} (a_n) \cap \partial_p
  \circ \cdots \partial_{n - 1} (b_n) \neq \emptyset
\end{equation*}
(not necessarily equal) then we glue along this part and defines a
composition
\begin{equation*}
  a_n \boks_p^n\; b_n.
\end{equation*}
See \cite{HOA}.
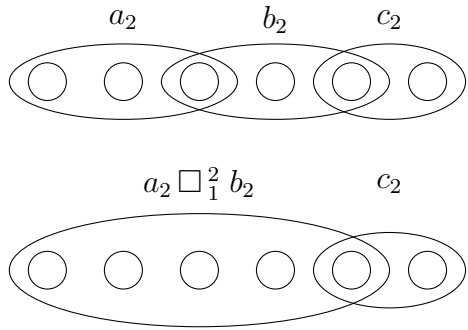
\begin{center}
  \begin{tikzpicture}
    \foreach \a in {0,1,2,3,4,5}{
      \draw (\a,0) circle (0.25cm);
    }

    \draw (1,0) ellipse (1.5cm and 0.5cm);
    \draw (3,0) ellipse (1.5cm and 0.5cm);
    \draw (4.5,0) ellipse (1cm and 0.5cm);

    \foreach \a/\alabel in {1/$a_2$,3/$b_2$,4.5/$c_2$}{
      \node at (\a,0.825){\alabel};
    }

    \begin{scope}[yshift=-2.5cm]
      \foreach \a in {0,1,2,3,4,5}{
        \draw (\a,0) circle (0.25cm);
      }

      \draw (2,0) ellipse (2.5cm and 0.75cm);
      \draw (4.5,0) ellipse (1cm and 0.5cm);

      \node at (2,1.15){$a_2 \boks_1^2\; b_2$};
      \node at (4.5,1.15){$c_2$};
    \end{scope}
  \end{tikzpicture}
\end{center}
Furthermore:
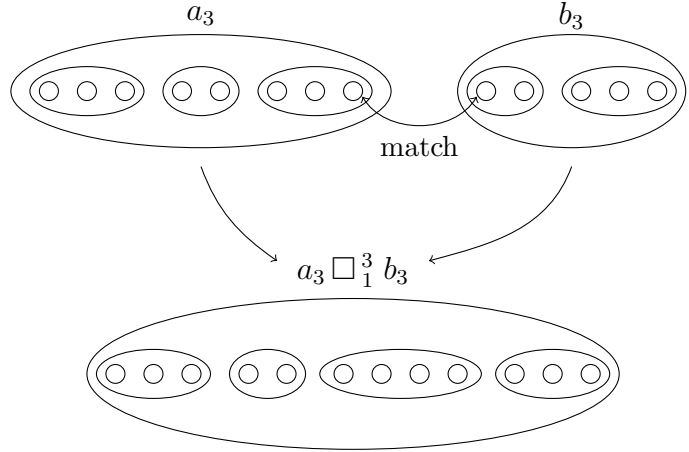
\begin{center}
  \begin{tikzpicture}
    \begin{scope}[scale=0.5]
      \foreach \a in {0,1,2,3.5,4.5,6,7,8}{
        \draw (\a,0) circle(0.25cm);
      }

      \draw (1,0) ellipse (1.5cm and 0.65cm);
      \draw (4,0) ellipse (1cm and 0.65cm);
      \draw (7,0) ellipse (1.5cm and 0.65cm);

      \draw (4,0) ellipse (5cm and 1.5cm);

      \node at (4,2){$a_3$};

      \foreach \b in {11.5,12.5,14,15,16}{
        \draw (\b,0) circle (0.25cm);
      }

      \draw (12,0) ellipse (1cm and 0.65cm);
      \draw (15,0) ellipse (1.5cm and 0.65cm);

      \draw (13.75,0) ellipse (3cm and 1.5cm);

      \node at (13.75,2){$b_3$};

      \draw[<->] (8.25,-0.15) to[out=-60,in=-120]
        node[midway,below,font=\small]{match} (11.25,-0.15);

      \draw[->] (4,-2) to[out=-70,in=145] (6,-4.5);
      \draw[->] (13.75,-2) to[out=-110,in=15] (10,-4.5);
    \end{scope}

    \begin{scope}[scale=0.5,xshift=1.75cm,yshift=-7.5cm]
      \foreach \c in {0,1,2,3.5,4.5,6,7,8,9,10.5,11.5,12.5}{
        \draw (\c,0) circle(0.25cm);
      }

      \draw (1,0) ellipse (1.5cm and 0.65cm);
      \draw (4,0) ellipse (1cm and 0.65cm);
      \draw (7.5,0) ellipse (2.125cm and 0.65cm);
      \draw (11.5,0) ellipse (1.5cm and 0.65cm);

      \draw (6.25,0) ellipse (7cm and 2cm);
      
      \node at (6.25,2.75){$a_3 \boks_1^3\; b_3$};
    \end{scope}
  \end{tikzpicture}
\end{center}
Fix level $i$ in the hyperstructure $\H(\C)$ and let $b_1,\ldots,b_k$
be bonds at level $i$, but gluable at level $j$, $j \leq i$.  Then we
form compositions
\begin{equation*}
  b_1 \boks_j^i\; b_2 \boks_j^i\; \cdots \boks_j^i\; b_k
\end{equation*}
which will form the $k$-simplexes at level $(i,j)$, let us call this
set $\Delta_k^{i,j}$.

We do this for all $i,j,k$ and form the smallest simplicial complex
having $\{\Delta_k^{i,j}\}$ as simplexes and we call this the nerve of
$\H(\C)$: 
\begin{equation*}
  N\H(\C).
\end{equation*}
This construction holds for any hyperstructure $\H$, not just $\H(\C)$
and extends the notion of classifying space or geometric realization
of $n$-categories.

In connection with the neural data it should be possible to compute
their homology and for suitable filtrations, their persistent
homology.  Our suggestion is that when we have data giving rise to a
code, we should look at the associated higher code or hyperstructure
$\H(\C)$.  For topological properties one should then study the nerve
$N\H(\C)$, not just $\Delta(\C)$.  This ought to capture the structure
of the data in a better way.

In a way forming the simplexes
\begin{equation*}
  \{b_1 \boks_j^i\; b_2 \boks_j^i\; \cdots \boks_j^i\; b_k\}
\end{equation*}
may amount to the same (or similar) process of ``geometric
realization'' of the sequence
\begin{equation*}
  \Delta(\C_k) \xleftarrow{\delta_{k - 1}^\ast} \Delta(\C_{k - 1})
  \leftarrow \cdots \xleftarrow{\delta_1^\ast} \Delta(\C_1)
\end{equation*}
making levelwise realization compatible with the $\delta_i^\ast$
maps.

Also in comparing two data sets this may be a useful approach.
Suppose that we have two data sets or codes, then we apply the
described procedure and produce two hyperstructures:
\begin{align*}
  \H &= \{\C_1,\ldots,\C_k\}\\
  \H' &= \{\C_1',\ldots,\C_k'\}
\end{align*}
with inbuilt boundary maps.  Then we may just compare bonds levelwise
by maps, are they for example one to one or what?

In general one may ask what kind of equivalence between $\H$ and $\H'$
is needed for $N(\H)$ and $N(\H')$ to have the same shape (or homotopy
type)?

It is tempting to think that neurons organized in a hierarchy of Hebb
assemblies may give rise to neural coding hyperstructures.  Such
assemblies may possibly be detected by the inference model for neural
data introduced in \cite{SDBB}, actually detecting tiers.

Hierarchies of Hebb assemblies, see \cite{Hebb,Scott} may
metaphorically compare to gene assemblies as in regulatory models of
the genome by Jacob and Monod \cite{Jacob} and Britten and Davidson
\cite{BD}.

Hyperstructures should play a role both in the analysis of neural data
and genomic data.  Let us give an illustration of how $N\H(\C)$ may
contain more information than $\Delta(\C)$.

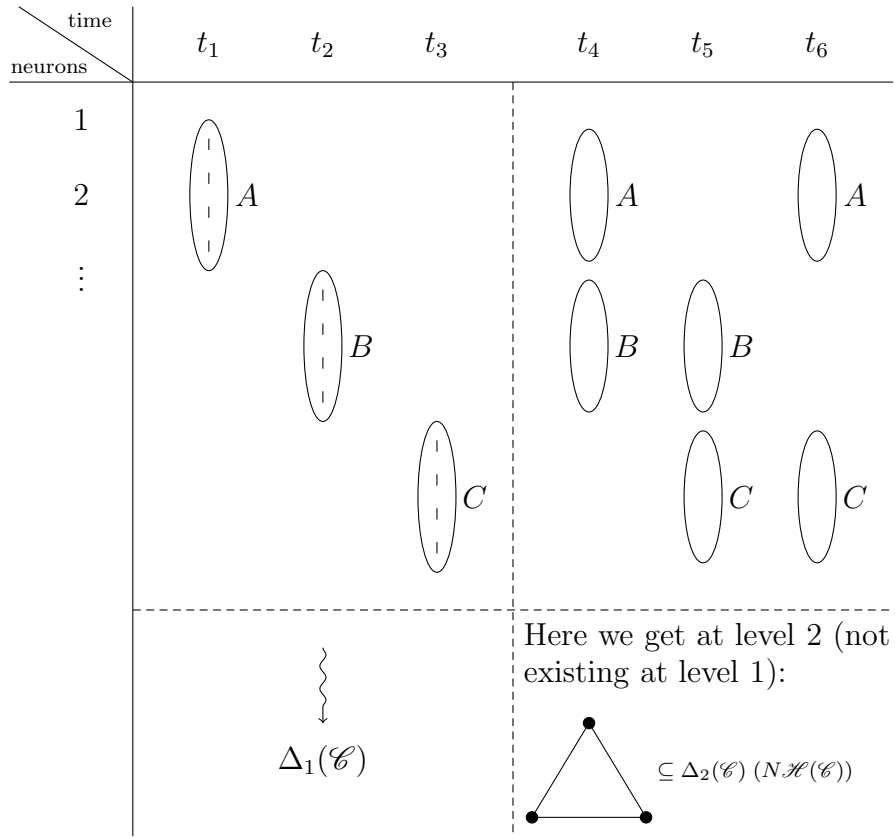
\begin{center}
  \begin{tikzpicture}
    \draw (-1.625,0) -- (10,0);
    \draw (0,1) -- (0,-10);
    \draw (0,0) -- (-1.5,1);

    \draw[densely dashed] (5,0) -- (5,-10);
    \draw[densely dashed] (0,-7) -- (10,-7);

    \node[right,font=\tiny] at (-1.75,0.175){neurons};
    \node[right,font=\tiny] at (-1,0.825){time};

    \foreach \t/\tlabel in {1/$t_1$, 2.5/$t_2$, 4/$t_3$, 6/$t_4$,
      7.5/$t_5$, 9/$t_6$}{
      \node at (\t,0.5){\tlabel};
    }

    \foreach \n/\nlabel in {-0.5/$1$,-1.5/$2$,-2.5/$\vdots$}{
      \node at (-0.675,\n){\nlabel};
    }

    \draw (1,-1.5) ellipse (0.25cm and 1cm);
    \foreach \A in {-0.75,-1.2,-1.65,-2.1}{
      \draw (1,\A) -- (1,\A-0.15);
    }
    \node at (1.5,-1.5){$A$};

    \draw (2.5,-3.5) ellipse (0.25cm and 1cm);
    \foreach \B in {-2.75,-3.2,-3.65,-4.1}{
      \draw (2.5,\B) -- (2.5,\B-0.15);
    }
    \node at (3,-3.5){$B$};

    \draw (4,-5.5) ellipse (0.25cm and 1cm);
    \foreach \C in {-4.75,-5.2,-5.65,-6.1}{
      \draw (4,\C) -- (4,\C-0.15);
    }
    \node at (4.5,-5.5){$C$};

    \draw (6,-1.5) ellipse (0.25cm and 0.875cm);
    \draw (6,-3.5) ellipse (0.25cm and 0.875cm);
    \node at (6.5,-1.5){$A$};
    \node at (6.5,-3.5){$B$};

    \draw (7.5,-3.5) ellipse (0.25cm and 0.875cm);
    \draw (7.5,-5.5) ellipse (0.25cm and 0.875cm);
    \node at (8,-3.5){$B$};
    \node at (8,-5.5){$C$};

    \draw (9,-1.5) ellipse (0.25cm and 0.875cm);
    \draw (9,-5.5) ellipse (0.25cm and 0.875cm);
    \node at (9.5,-1.5){$A$};
    \node at (9.5,-5.5){$C$};

    \draw[->,decorate,decoration={snake,amplitude=.4mm,segment
        length=3mm,post length=1mm}] (2.5,-7.5) -- (2.5,-8.5);
    \node at (2.5,-9){$\Delta_1(\C)$};

    \node[below right,text width=5cm] at (5,-7){Here we get at level
      $2$ (not existing at level $1$):};
    \draw (5.25,-9.75) -- (6,-8.5) -- (6.75,-9.75) -- cycle;
    \filldraw[fill=black] (5.25,-9.75) circle(0.075cm);
    \filldraw[fill=black] (6,-8.5) circle(0.075cm);
    \filldraw[fill=black] (6.75,-9.75) circle(0.075cm);
    \node[right,font=\tiny] at (6.75,-9.15){$\subseteq \Delta_2(\C) \;
      (N\H(\C))$};
  \end{tikzpicture}
\end{center}
In addition come lower level gluing simplexes. For one level, see \cite{EK}.

\subsection*{Acknowledgements}
I would like to thank M.~Thaule for his kind technical assistance in
preparing the manuscript.

\end{document}